\documentclass[12pt,a4paper]{article}
\usepackage{amsmath,amssymb,amsthm}
\usepackage[dvips,final]{epsfig}
\usepackage[dvips,final]{graphicx}
\usepackage{color}
\usepackage{soul}

\usepackage{xcolor} 

\definecolor{grayl}{rgb}{0.94,0.94,0.94}
\definecolor{gre}{rgb}{0.2,0.6,0}
\definecolor{vio}{rgb}{.6,0,0.6}
\definecolor{ora}{rgb}{1,.5,0}
\definecolor{cya}{rgb}{0,.8,.8}

\def\grayl#1{\textcolor{grayl}{#1}}

\newtheorem{theorem}{Theorem}[section]
\newtheorem{hyp}[theorem]{Hypothesis}
\newtheorem{lemma}[theorem]{Lemma}
\newtheorem{corollary}[theorem]{Corollary}
\newtheorem{proposition}[theorem]{Proposition}
\theoremstyle{definition}
\newtheorem{definition}[theorem]{Definition}
\newtheorem{remark}[theorem]{Remark}

\numberwithin{equation}{section}

\usepackage{titlesec}
\titleformat{\section}{\bfseries}{\thesection}{1em}{}
\titleformat{\subsection}{\itshape}{\thesubsection}{1em}{}

\newcommand{\real}{\mathbb{R}}
\newcommand{\nat}{\mathbb{N}}

\newcommand{\Var}{\mathop{\mathrm{Var}}}
\def\expe{\mathrm{e}}
\def\dd{\,\mathrm{d}}
\def\dist{\mathrm{dist}}
\def\meas{\mathrm{meas}}
\def\Int{\mathrm{Int\,}}

\def\duw#1{\left\langle\!\left\langle #1 \right\rangle\!\right\rangle}
\def\scal#1{\left\langle #1 \right\rangle}
\def\ve{\varepsilon}
\def\ap{\alpha}
\def\L{\mathcal{L}}

\def\for{\mbox{ for }\ }

\def\be{\begin{equation}\label}
\def\ee{\end{equation}}

\def\doxi{\dot \xi}
\def\doxx{\dot \xi_2}
\def\dota{\dot \eta}
\def\dou{\dot u}

\def\dow{\dot w}
\def\dox{\dot x}

\def\barr{\begin{array}}
\def\earr{\end{array}}

\def\bearr{\begin{eqnarray}}
\def\eearr{\end{eqnarray}}

\def\bears{\begin{eqnarray*}}
\def\eears{\end{eqnarray*}}

\newcommand{\AC}{W^{1,1}(0,T; X)}
\newcommand{\domab}{\AC \times W^{1,1}(0,T\,; W)}

\newfont{\ctv}{msam10}

\newcommand{\bbox}{\mbox{\ctv \symbol{4}}}
\def\QED{{${}\hfill\bbox$}}
\newenvironment{pf}[1]{\par\vskip1mm{\noindent\it #1.}\ }{\QED\par
\vskip2mm}

\def\bpf{\begin{pf}}
\def\epf{\end{pf}}

\font\prox = mfprox1
\font\bfi = cmbxti10 scaled \magstep2

\newenvironment{pkrev}{\color{vio}}{\color{black}}
\newcommand{\bpk}{\begin{pkrev}}
\newcommand{\epk}{\end{pkrev}}

\newenvironment{pkrevg}{\color{gre}}{\color{black}}
\newcommand{\bpg}{\begin{pkrevg}}
\newcommand{\epg}{\end{pkrevg}}

\begin{document}

\title{Explicit and implicit non-convex sweeping processes in the space of absolutely continuous functions\footnote{Supported by the GA\v CR Grant No.~20-14736S, RVO: 67985840, and by the European Regional Development Fund, Project No. CZ.02.1.01/0.0/0.0/16{\_}019/0000778.}}

\date{}

\author{Pavel Krej\v c\'\i\footnote{Faculty of Civil Engineering, Czech Technical University, Th\'akurova 7, 16629 Praha 6, Czech Republic (\texttt{pavel.krejci@cvut.cz})}, $\;$ 
Giselle Antunes Monteiro\footnote{Institute of Mathematics, Czech Academy of Sciences, \v Zitn\'a 25,
11567~Praha 1, Czech Republic,  (\texttt{gam@math.cas.cz})}, $\;$ 
Vincenzo Recupero\footnote{Dipartimento di Scienze Matematiche, Politecnico di Torino, Corso Duca degli Abruzzi 24, 10129 Torino, Italy (\texttt{vincenzo.recupero@polito.it})} \footnote{Vincenzo Recupero is a member of GNAMPA-INdAM.}}


\maketitle

\begin{abstract}
We show that sweeping processes with possibly non-convex prox-regular constraints generate a strongly continuous input-output mapping in the space of absolutely continuous functions. Under additional smoothness assumptions on the constraint we prove the local Lipschitz continuity of the input-output mapping. Using the Banach contraction principle, we subsequently prove that also the solution mapping associated with the state-dependent problem is locally Lipschitz continuous.

\ 

\begin{footnotesize}
\end{footnotesize}


\end{abstract}

\section*{Introduction}

The present paper is a continuation of \cite{kmr}, where we have studied a class of constrained evolution problems, called \emph{sweeping processes}, in the framework of a real Hilbert space $X$ endowed with scalar product $\scal{x,y}$ and norm $|x| = \sqrt{\scal{x,x}}$. In order to describe the processes studied in \cite{kmr}, we assume that we are given right-continuous functions $u:[0,T] \to X$ and $w:[0,T] \to W$, where $W$ is a real Banach space and we suppose that $u$ and $w$ are regulated, i.\,e. they admit left limits at every point $t \in (0,T]$. We also assume that a not necessarily convex moving constraint $Z(w(t))\subset X$ is given and that $Z(w(t))$ is $r$-prox-regular, i.e. $Z(w(t))$ is a closed sets having a neighborhood of radius $r > 0$ where the metric projection exists and is unique.

Assuming that the sets $Z(w(t))$ satisfy a suitable uniform non-empty interior condition, 
we have proved in \cite{kmr} that
for every initial condition $x_0 \in Z(w(0))$ there exists a right-continuous function $\xi:[0,T] \to X$ of bounded variation ($BV$) such that the variational inequality
\be{s4}
\int_0^T \scal{x(t) - z(t), \dd\xi(t)} + \frac1{2r}\int_0^T |x(t) - z(t)|^2 \dd V(\xi)(t) \ge 0, \quad x(0) = x_0,
\ee
is satisfied for every regulated test function $z:[0,T] \to X$ such that $z(t) \in Z(w(t))$ for all $t \in [0,T]$, with $x(t) := u(t) - \xi(t)$ and $V(\xi)(t) := \Var_{[0,t]} \xi$, the variation of $\xi$ over $[0,t]$ for $t \in [0,T]$.
The two integrals in \eqref{s4} can be interpreted 
in the sense of the Kurzweil integral introduced in \cite{kur57}: In the first integral we are integrating $X$-valued functions, while the second integral corresponds to the standard case of real-valued functions.

Since the normal cone $N_Z(x)$ of a closed $Z \subseteq H$ at $x \in Z$ is defined by the formula
\be{e2}
N_Z(x) = \left\{\xi \in X\ :\  \scal{\xi, x-z} + \frac{|\xi|}{2r} |x-z|^2 \ge 0 \quad \forall z \in Z\right\},
\ee
the variational inequality \eqref{s4} can be formally interpreted as a $BV$ integral formulation of the differential inclusion
\begin{equation}
  \dot \xi(t) \in - N_{C(t)}(\xi(t)), \quad \xi(0) = u(0) - x_0
\end{equation}
with $C(t) = u(t) - Z(w(t))$, $t \in [0,T]$.

In \cite[Section 5]{kmr}, we have shown under some technical assumptions, 
but dropping the uniform non-empty interior condition for $Z(w(t))$, that if the inputs $u$, $w$ are absolutely continuous, then the output $\xi$ is absolutely continuous and satisfies the pointwise variational inequality
\be{a1}
\scal{x(t) - z, \dot\xi(t)} + \frac{|\dot\xi(t)|}{2r}|x(t) - z|^2 \ge 0, \quad x(t) + \xi(t) = u(t), \quad x(0) = x_0
\ee
for a.\,e. $t \in (0,T)$ and all $z \in Z(w(t))$. The existence and uniqueness result for \eqref{a1} was stated and proved in \cite[Corollary 5.3]{kmr} and we recall the precise statement below in Proposition \ref{ap1}.

A detailed survey of the literature related to non-convex sweeping processes was given in \cite{kmr} and we do not repeat it here. Instead, we pursue further the study of \eqref{a1} in the space of absolutely continuous functions. Let us mention only the publications that have particularly motivated our research, namely the pioneering paper \cite{moreau} where the concept of sweeping process was elaborated, the detailed 
studies \cite{csw95, prt} of prox-regular sets, and a deep investigation of prox-regular sweeping processes carried out in \cite{cmm,Thi03}.

It turns out that it is convenient in this context to represent the sets $Z(w) = \{x \in X: G(x,w) \le 1\}$ as sublevel sets of a function $G: X \times W \to [0,\infty)$ satisfying suitable technical assumptions. 
A detailed comparison of different continuity criteria has been done in the convex case in \cite{bks}. In the nonconvex case treated in the present paper we prove as our main result that the input-output mapping $(u,w) \mapsto \xi$
is strongly continuous with respect to the $W^{1,1}$-norms if $G$ is continuously differentiable with respect to both $x$ and $w$, and Lipschitz continuous if both gradients $\nabla_x G, \nabla_w G$ are Lipschitz continuous. 
As a consequence of the Lipschitz input-output dependence, we apply the Banach contraction principle to prove the unique solvability of an implicit state dependent problem with $w$ of the form $w(t) = g(t, u(t), \xi(t))$ with a given smooth function $g:[0,T]\times X \times X \to W$.
The authors are not aware of any result of this kind in the literature on prox-regular sweeping processes. Implicit problems in the convex case have been solved under suitable additional compactness assumptions in \cite{kmm97, kmm98} and without compactness in \cite{bks}. The non-convex case has been considered for example in \cite{AdlHad18, AdlHadLe19, jv, jv19, msz}, but to our knowledge, in all existing publications, the sweeping process is regularized by some kind of compactification or viscous regularization. In our case, no compactification or other kind of regularization comes into play.

The paper is structured as follows. In Section \ref{pro}, we identify sufficient conditions on the function $G(\cdot, w)$ which guarantee that the sublevel set $Z(w)$ is $r$-prox-regular for every $w \in W$. Section \ref{abs} is devoted to finding additional hypotheses
on the $w$-dependence of $G$ which guarantee the validity of the existence and uniqueness result for Problem \eqref{a1} in \cite[Corollary 5.3]{kmr}. The strong continuity of the $(u,w) \mapsto \xi$ input-output mapping with respect to the $W^{1,1}$-norm is proved in Section \ref{lip}, and the local Lipschitz continuity of the mapping $u \mapsto \xi$ in the implicit case $w(t) = g(t, u(t), \xi(t))$ is proved in Section \ref{imp}.


\section{Prox-regular sets of class $C^1$}\label{pro}

Let us start by recalling the definition of prox-regular set, in agreement with 
\cite[items (a) and (g) of Theorem 4.1]{prt}.

\begin{definition}\label{d1}
Let $X$ be a real Hilbert space endowed with scalar product $\scal{\cdot, \cdot}$ and norm 
$|x| = \sqrt{\scal{x,x}}$, let $Z\subset X$ be a closed connected set, and let $\dist(x,Z) := \inf\{|x-z|: z\in Z\}$ denote the distance of a point $x\in X$ from the set $Z$. Let $r>0$ be given. We say that $Z$ is 
{\em $r$-prox-regular} if the following condition hold.
\be{e1}
\forall y \in X:\ \dist(y,Z) = d \in(0, r) \ \ \exists\, x \in Z: \dist\left(x+\frac{r}{d}(y-x),Z\right) = \frac{r}{d}|y-x|= r.
\ee
\end{definition}

We have the following characterization of prox-regular sets (see, e.\,g., \cite{prt,kmr}).

\begin{lemma}\label{l2}
A set $Z\subset X$ is $r$-prox-regular if and only if for every $y \in X$ such that $d =\dist(y,Z) < r$ there exists a unique $x\in Z$ such that $|y-x| = d$ and
\be{e2w}
\scal{y-x, x-z} + \frac{|y-x|}{2r} |x-z|^2 \ge 0 \quad \forall z \in Z.
\ee
\end{lemma}

We represent the sets $Z$ as the sublevel sets of a function $G:X \to [0,\infty)$ in the form
\be{w1}
Z = \{z \in X: G(z) \le 1\}\,.
\ee
We define the gradient $\nabla G(x) \in X$ of $G$ at a point $x \in X$ by the formula
\be{gg}
\scal{\nabla G(x), y} = \lim_{t\to 0}\frac1t\big(G(x + ty) - G(x)\big)\quad \forall y \in X\,.
\ee
The following hypothesis is assumed to hold.

\begin{hyp}\label{clh}
Let $X$ be a real Hilbert space endowed with scalar product $\scal{\cdot, \cdot}$ and norm 
$|x| = \sqrt{\scal{x,x}}$. We assume that \eqref{w1} holds for a function $G:X \to [0,\infty)$, $\nabla G(x)$ exists for every $x \in X$, and there exist positive constants $\lambda, c$ and a continuous increasing function $\mu: [0,\infty) \to [0,\infty)$ such that $\mu(0) = 0$, $\lim_{s \to \infty} \mu(s) = \infty$, and
\begin{itemize}
\item[{\rm (i)}] $G(x) = 1\ \Longrightarrow\ |\nabla G(x)| \ge c>0$ for all $x \in X$;
\item[{\rm (ii)}] $|\nabla G(x) - \nabla G(y)| \le \mu(|x-y|)$ for all $x,y \in Z$; 
\item[{\rm (iii)}] 
$\scal{\nabla G(x) - \nabla G(z), x-z} \ge -\lambda |x-z|^2$ \ for all $x\in \partial Z$ and $z \in Z$.
\end{itemize}
\end{hyp}

Throughout the paper, for a set $S \subset X$, the symbols $\partial S$, $\Int S$, and $\overline{S}$ will denote respectively the boundary, the interior, and the closure of $S$. It is easy to check that under Hypothesis \ref{clh} we have $\partial Z = \{x \in X: G(x) = 1\}$, so that an element $x \in X$ belongs to 
$\Int Z$ if and only if $G(x) < 1$. Indeed, if $G(x) = 1$, then
$$
\lim_{t\to 0} \frac1{t} \big(G(x + t \nabla G(x)) - G(x)\big) = |\nabla G(x)|^2 \ge c^2 > 0.
$$
For $n\in \nat$ put $x_n := x + \frac1n\nabla G(x)$. We have $G(x_n) > 1$ for $n$ sufficiently large, hence $x_n \notin Z$. Since $x_n$ converge to $x$ as $n\to \infty$, we conclude that $x \in \partial Z$.

In the convex case, we can choose $G$ to be the Minkowski functional $M_Z$ (or gauge) associated with $Z$ defined as $M_Z(x) = \inf\{s>0: \frac1s x \in Z\}$. Then condition (iii) of Hypothesis \ref{clh} is automatically satisfied, since $\nabla M_Z$ is monotone, and (ii) is just the uniform continuity condition of $\nabla M_Z$. For non-convex sets $Z$, condition (iii) excludes sharp concavities of $\partial Z$.

We now prove the following result.

\begin{proposition}\label{lp1}
Let Hypothesis \ref{clh} hold and let $r = c/\lambda$.
Then for every $y \in X$ such that $d =\dist(y,Z) \in (0,r)$ there exists a unique $x\in \partial Z$ such that 
$$
y = x+ d\frac{\nabla G(x)}{|\nabla G(x)|},
$$ 
and
\be{e2wg}
\scal{\frac{\nabla G(x)}{|\nabla G(x)|}, x-z} + \frac{\lambda}{2c} |x-z|^2 \ge 0 \quad \forall z \in Z.
\ee
In particular, $Z$ is $r$-prox-regular.
\end{proposition}

The statement of Proposition \ref{lp1} is not new. The finite-dimensional case was already solved in \cite{Via83}.
The fact that the conditions of Hypothesis 1.3 are sufficient for a set given
by (1.3) to be prox-regular also in the infinite-dimensional case was shown in \cite[Theorem 9.1]{AdlHadThi17} (see also \cite{AdlHadThi16}).
The proof there refers to a number of deep concepts from non-smooth analysis along the lines, e.\,g., of \cite[Chapter 2]{clsw98}.
Here we present instead an elementary self-contained proof using no other analytical tools but the properties of the scalar product, and the argument is split into several steps including two auxiliary Lemmas. 

\begin{lemma}\label{lw1}
Let Hypothesis \ref{clh} hold and let $r = c/\lambda$. Then for all $x \in \partial Z$ and $z\in Z$ we have
$$
\scal{\nabla G(x),x-z} + \frac{|\nabla G(x)|}{2r} |x-z|^2 \ge 0.
$$
\end{lemma}

\bpf{Proof of Lemma \ref{lw1}}
For $x \in \partial Z$ and $z \in Z$ we have 
\begin{align*}
0 & \le G(x) - G(z) = \frac{\dd}{\dd t} \int _0^1 G(z + t(x-z))\dd t = \int_0^1 \scal{\nabla G(z + t(x-z)),x-z}\dd t\\
&= \scal{\nabla G(x),x-z}
- \int_0^1 \scal{\nabla G(x) - \nabla G(x - (1{-}t)(x-z)),x-z}\dd t\\
&\le \scal{\nabla G(x),x-z} + \frac{\lambda}{2} |x-z|^2.
\end{align*}
and the assertion follows.
\epf

\begin{lemma}\label{cl}
Let Hypothesis \ref{clh} hold, let $V \subset X$ be the set of all $y \in X$ for which there exists $x \in Z$ such that $|y-x| = \dist(y,Z)$, and let $U_r := \{y \in Z: \dist(y, Z) < r\}$. Then the set $V\cap U_r$ is dense in 
$\overline{U_r}$.
\end{lemma}

A highly involved proof of Lemma \ref{cl} can be found in a much more general setting, e.\,g., in 
\cite[Theorem 3.1, p. 39]{clsw98}. For the reader's convenience, we show that in our special case, it can be proved in an elementary way.

\begin{figure}[htb] \label{f1}
\begin{center}
\begin{picture}(75,50)
\put(0,0){\grayl{\prox \char6}}
\put(0,0){\prox \char5}
\put(32,24.5){\small $x_\ap$}
\put(12,33){\small $x_\ap + t(y-y_\ap)$}
\put(76,34){\small $y$}
\put(76,14){\small $y_\ap$}
\put(14,10){\bfi Z}
\end{picture}
\caption{Illustration to Lemma \ref{cl}.}
\end{center}
\end{figure}
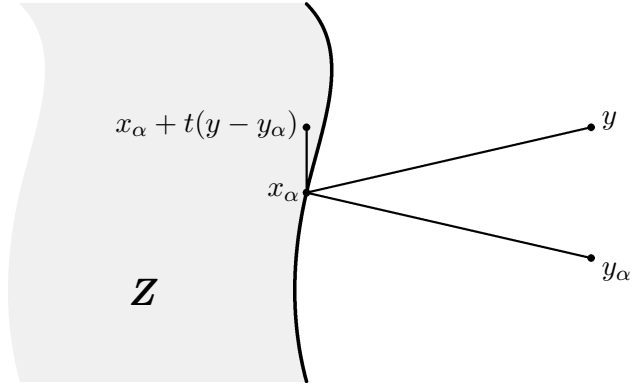

\bpf{Proof of Lemma \ref{cl}}
We prove that for a given $y \in U_r$ and every $\ve>0$ there exists $y^* \in V\cap U_r$ such that
\be{cl0}
|y - y^*| <\ve.
\ee
Let $y \in X$ be arbitrarily chosen such that $d := \dist(y,Z) < r$. For any $\ap \in (0,r-d)$ we find  
$x_\ap \in \partial Z$ such that $|y-x_\ap| = d+\ap$ and put (see Figure \ref{f1})
\be{cl1}
y_\ap := x_\ap + (d+\ap)n_\ap, \ \ \gamma_\ap := |\nabla G(x_\ap)|, \ \ 
n_\ap := \frac{\nabla G(x_\ap)}{\gamma_\ap}, \ \ \bar n_\ap := \frac{y-x_\ap}{d+\ap}, \ \  
q_\ap := \scal{n_\ap, \bar n_\ap}.
\ee


Using Lemma \ref{lw1} we check that $y_\ap \in V$, $|y_\ap - x_\ap| = d+\ap = \dist(y_\ap, Z)$. We have by \eqref{cl1} that $y - y_\ap = (d+\ap)(\bar n_\ap - n_\ap)$, hence,
\begin{equation}\label{cl1-2}
  |y - y_\ap|^2 = (d+\ap)^2(2 - 2\scal{\bar n_\ap, n_\ap}) = 2(d+\ap)^2(1 - q_\ap).
\end{equation}
This implies that $-1 \le q_\ap < 1$. Indeed, if $q_\ap =1$, then $y = y_\ap$ and $\dist(y, Z) = d+\ap$ which is a contradiction. Furthermore,
\begin{align}\nonumber
\lim_{t\to 0+} \frac1t\big(G(x_\ap + t(y-y_\ap)) - G(x_\ap)\big) 
&= \scal{\nabla G(x_\ap),y-y_\ap}  =
\gamma_\ap (d+\ap)\scal{n_\ap, \bar n_\ap - n_\ap} \\ \label{cl2}
&= \gamma_\ap (d+\ap)(q_\ap - 1) < 0.
\end{align}
We have $G(x_\ap) = 1$, hence $G(x_\ap + t(y-y_\ap)) < 1$ in a right neighborhood of $0$. Put 
$T_\ap := \inf\{t > 0\, : G(x_\ap + t(y-y_\ap)) \ge 1\}$. Then $T_\ap > 0$, and for every $t_\ap\in (0,T_\ap)$, by hypothesis that $\dist(y, Z) = d$ and by \eqref{cl1}, \eqref{cl1-2} we find 
\begin{align*}
d^2 &< |y -(x_\ap + t_\ap(y-y_\ap))|^2 = |(1-t_\ap)(y-x_\ap) + t_\ap(y_\ap - x_\ap)|^2\\
&= (d+\ap)^2|(1-t_\ap)\bar n_\ap + t_\ap n_\ap|^2 \\
&= (d+\ap)^2\big((1-t_\ap)^2 + t_\ap^2 + 2t_\ap(1-t_\ap) q_\ap\big)\\
&= (d+\ap)^2\big(1 -2t_\ap(1-t_\ap)(1- q_\ap)\big) \\
& \ = (d+\ap)^2 - t_\ap(1-t_\ap)|y-y_\ap|^2, 
\end{align*}
so that
\begin{equation}\label{c12-2}
  t_\ap(1-t_\ap)|y-y_\ap|^2 < (d+\ap)^2 - d^2 = (2d+\alpha)\alpha \quad \forall t_\ap \in (0,T_\ap), 
    \ \forall \alpha \in (0,r-d).
\end{equation}
If $\limsup_{\alpha \searrow 0} T_\ap > 1/2$, then for all $\alpha$ such that $T_\ap > 1/2$  we can take $t_\ap=1/2$ in \eqref{c12-2} and obtain that
\begin{align*}
|y - y_\ap|^2  <  4\ap(2d+\ap), 
\end{align*}
and \eqref{cl0} is satisfied provided we choose $y^* = y_\ap$ for a sufficiently small $\ap$ such that
$T_\ap > 1/2$.

It remains to consider the case $\limsup_{\alpha \searrow 0} T_\ap \le 1/2$. We have 
$G(x_\ap + T_\ap(y-y_\ap)) = 1$, thus, by Hypothesis \ref{clh} and \eqref{cl1-2}, we find
\begin{align*}
0 &= G(x_\ap + T_\ap(y-y_\ap)) - G(x_\ap) = \int_0^{T_\ap}\scal{\nabla G(x_\ap + t(y-y_\ap)), y-y_\ap}\dd t\\
&= T_\ap \scal{\nabla G(x_\ap),y-y_\ap} + \int_0^{T_\ap}\scal{\nabla G(x_\ap + t(y-y_\ap)) - \nabla G(x_\ap), y-y_\ap}\dd t\\
&\le -T_\ap \gamma_\ap (d+\ap)(1-q_\ap) + \int_0^{T_\ap} 
\mu(t|y-y_\ap|)|y-y_\ap|\dd t\\
&= -\frac{T_\ap \gamma_\ap|y-y_\ap|^2}{2(d+\ap)} + \int_0^{T_\ap|y-y_\ap|}
\mu(\sigma)\dd \sigma.
\end{align*}
For $p \ge 0$ put
$$
M(p) = \int_0^{p} \mu(\sigma)\dd \sigma.
$$
The function $M:[0,\infty) \to [0,\infty)$ is increasing and convex, and the function $\hat M(p) = M(p)/p$ is increasing, unbounded, and $\hat M(0+) = \mu(0) = 0$. From the above computations we conclude that
\be{cl3}
\hat M(T_\ap|y-y_\ap|) \ge \frac{c}{2r}|y-y_\ap|,
\ee
where $c$ is the constant from Hypothesis \ref{clh}\,(i).

Put $t_\ap = T_\ap/2$. Then $t_\ap < 1/2$ for $\ap$ sufficiently small. Using \eqref{cl3} and \eqref{c12-2} we infer that
\begin{align}
|y-y_\ap|\,\hat M^{-1} \left(\frac{c}{2r}|y-y_\ap|\right) & \le T_\ap |y-y_\ap|^2 = 2t_\ap |y-y_\ap|^2 
\notag \\
& < 2(1-t_\ap)^{-1}(2d+\ap)\ap < 4(2d+\ap)\ap,
\end{align}
hence, putting $p_\ap := \hat M^{-1} (c/(2r) |y-y_\ap|)$, we obtain
$$
 M(p_\ap) = p_\ap \hat M(p_\ap) < \frac{2c}{r}(2d+\ap)\,\alpha
$$
and we obtain \eqref{cl0} for $y^* = y_\ap$ and $\ap>0$ sufficiently small from the continuity of $\mu$, $M$ and $\hat{M}$ at $0$.
\epf

We are now ready to prove Proposition \ref{lp1}.

\bpf{Proof of Proposition \ref{lp1}}
To prove that $Z$ is $r$-prox-regular, consider any $d \in (0,r)$ and put
\begin{align*}
\Gamma &= \partial Z =\{y \in X: G(y) = 1\},\\
\Gamma_d &= \{y \in X: \dist(y,Z) = d\}.
\end{align*}
Let $f: \Gamma \to X$ be the mapping defined by the formula
\be{pra}
f(x) = x + d\frac{\nabla G(x)}{|\nabla G(x)|}.
\ee
Then $f$ is continuous and $f(\Gamma) \subset \Gamma_d$. Indeed, we have $|f(x) - x| = d$ and, choosing an arbitrary $z \in Z$,
$$
|f(x) - z|^2 = |x-z|^2 + d^2 + 2d\scal{\frac{\nabla G(x)}{|\nabla G(x)|}, x-z} \ge \left(1 - \frac{d}{r}\right)|x-z|^2 + d^2 \ge d^2
$$
by virtue of Lemma \ref{lw1}. Furthermore, by Lemma \ref{cl}, for $y$ from a dense subset of $\Gamma_d$ there exists $x \in \Gamma$ and a unit vector $n(x)$ such that $y= x+ dn(x)$ and $|y-z| \ge d$ for all $z \in Z$. We consequently have $|x+tn(x) - z| \ge |x+dn(x) - z| - (d-t) \ge t$ for all $t \in (0,d]$. Put
$$
\hat n(x) = \frac{\nabla G(x)}{|\nabla G(x)|}.
$$
We have for all $t \in (0,d)$ that $x + tn(x) - t\hat n(x) \notin \Int Z$, hence 
$G(x + tn(x) - t\hat n(x)) \ge 1$, 
and
$$
0 \le \lim_{t\to 0+} \frac1t\left(G(x + tn(x) - t\hat n(x)) - G(x)\right) = \scal{\nabla G(x), n(x) -\hat n(x)},
$$
and we conclude that $n(x) =\hat n(x)$. The range $f(\Gamma)$  of the mapping $f$ defined by \eqref{pra} is therefore dense in $\Gamma_d$. Assume that there exists $y \in \Gamma_d\setminus f(\Gamma)$. We find a sequence of elements $y_j \in f(\Gamma)$, $j \in \nat$, which converges to $y$, $y_j = f(x_j)$. For $j,k \in \nat$ we have in particular
$$
y_j - y_k = x_j - x_k + d \left(\frac{\nabla G(x_j)}{|\nabla G(x_j)|}-\frac{\nabla G(x_k)}{|\nabla G(x_k)|}\right).
$$
By Lemma \ref{lw1} we have
$$
\scal{\left(\frac{\nabla G(x_j)}{|\nabla G(x_j)|}-\frac{\nabla G(x_k)}{|\nabla G(x_k)|}\right), x_j - x_k} 
\ge -\frac{1}{r} |x_j - x_k|^2,
$$
hence,
$$
|x_j - x_k|^2 \le \frac{r}{r-d} \scal{y_j - y_k,x_j - x_k}.
$$
We conclude that $\{x_j\}$ is a Cauchy sequence in $X$, hence it converges to some $x \in \Gamma$ and the continuity of $f$ yields $y = f(x)$. We have thus proved that for each $y \in \Gamma_d$ there exists a unique $x \in \Gamma$ such that $y = f(x)$, and the assertion follows from Lemma \ref{cl}.
\epf


\section{Absolutely continuous inputs}\label{abs}

We now consider a family of sets $\{Z(w): w\in W\}$ parameterized by elements $w$ of a Banach space $W$ 
with norm $|\cdot|_W$ and defined as the sublevel sets
\begin{equation}\label{zw}
  Z(w) = \{z \in X: G(x,w) \le 1\}
\end{equation}
of a locally Lipschitz continuous function $G: X \times W \to [0,\infty)$. Similarly as in \eqref{gg}, we define the partial gradients $\nabla_x G(x,w) \in X$, $\nabla_w G(x,w) \in W'$ for $x\in X$ and $w \in W$ by the identities
\begin{align}\label{ggx}
\scal{\nabla_x G(x,w), y} &= \lim_{t\to 0}\frac1t\big(G(x+ty,w) - G(x,w)\big)\quad \forall y \in X\,,\\ \label{ggw}
\duw{\nabla_w G(x,w), v} &= \lim_{t\to 0}\frac1t\big(G(x,w+tv) - G(x,w)\big)\quad \forall v \in W\,,
\end{align}
where $W'$ is the dual of $W$, and $\duw{\cdot,\cdot}$ is the duality $W \to W'$.
We assume the following hypothesis to hold.

\begin{hyp}\label{hl}
Let $X$ be a real Hilbert space endowed with scalar product $\scal{\cdot, \cdot}$ and norm 
$|x| = \sqrt{\scal{x,x}}$ and let $W$ be a real Banach space with norm $|\cdot|_W$. We assume that 
\eqref{zw} holds for a locally Lipschitz continuous function $G: X \times W \to [0,\infty)$ for which 
$\nabla_x G(z,w)$ exists for every $(z,w) \in Z\times W$ and there exist positive constants $\lambda, c, L$ and  
functions $\mu_1: W \times[0,\infty) \to [0,\infty)$, $\mu_2: [0,\infty) \to [0,\infty)$ such that 
$\mu_1(w,0) = \mu_2(0) = 0$, $\lim_{s \to \infty} \mu_1(w,s) = \lim_{s \to \infty} \mu_2(s) = \infty$
for every $w \in W$, and 
\begin{itemize}
\item[{\rm (i)}] $G(x,w) = 1\ \Longrightarrow\ |\nabla_x G(x,w)| \ge c>0$ for all $x \in X$ and $w \in W$;
\item[{\rm (ii)}] $|\nabla_x G(x,w) - \nabla_x G(y,w)| \le \mu_1(w,|x-y|)$ for all $x,y \in Z(w)$ and $w \in W$;
\item[{\rm (iii)}] 
$\scal{\nabla_x G(x,w) - \nabla_x G(z,w), x-z} \ge -\lambda |x-z|^2$ \ for all $x\in \partial Z(w)$, 
$z \in Z(w)$, and $w \in W$;
\item[{\rm (iv)}] $|G(x,w) - G(x,w')| \le L |w - w'|_W$ for all $x \in X$ and $w,w' \in W$;
\item[{\rm (v)}] $\forall \rho >0 \ \forall w \in W \ \forall x \in X:$
\be{da}
\dist(x, Z(w)) \ge \rho\ \Longrightarrow\ G(x,w)-1 \ge \mu_2(\rho).
\ee
\end{itemize}
\end{hyp}

The property (v) in Hypothesis \ref{hl} is a kind of uniform coercivity of the function $G$ which will play a role in the next Lemma.
Let us observe that Hypothesis \ref{hl} and Proposition \ref{lp1} imply that the set $Z(w)$ given by \eqref{zw} is $r$-prox-regular for every $w \in W$.

\begin{lemma}\label{ldh}
Let Hypothesis \ref{hl} hold. Then for every $K>0$ there exists a constant $C_K>0$ such that
\be{lipw}
\max\{|w_1|_W, |w_2|_W\} \le K \ \Longrightarrow\  d_H(Z(w_1),Z(w_2)) \le  C_K |w_1 - w_2|_W, 
\ee
for every $w_1,w_2 \in W$, where $d_H$ denotes the Hausdorff distance
$$
d_H(Z(w_1),Z(w_2)):= \max\Big\{\sup_{z \in Z(w_1)}\dist(z, Z(w_2)), \sup_{z' \in Z(w_2)}\dist(z', Z(w_1))\Big\}.
$$
\end{lemma}

\bpf{Proof}
Let $K>0$ be given and let $\max\{|w_1|_W, |w_2|_W\} \le K$. We first check that
\be{ek}
d_H(Z(w_1),Z(w_2)) \le D_K:= \mu_2^{-1}(2KL).
\ee
Indeed, if this was not true, we can assume that there exists $x \in Z(w_1)$ such that  
$\dist(x, Z(w_2)) \ge D_K + \ap$ for some $\ap > 0$. By Hypotheses \ref{hl}\,(iv)-(v) we then have
$$
L |w_1 - w_2|_W \ge G(x, w_2) - G(x,w_1) \ge G(x, w_2) - 1 \ge \mu_2(D_K + \ap) > 2KL, 
$$
which is a contradiction. 

We now consider the cases $d_H(Z(w_1),Z(w_2))\ge r$ or $d_H(Z(w_1),Z(w_2)) < r$ separately. Let us start with the case
\\[2mm]
A. $d_H(Z(w_1),Z(w_2)) \ge r$.
\\[2mm]
Then for $x \in Z(w_1)$ we have by Hypotheses \ref{hl}\,(iv)-(v) that
$$
L |w_1 - w_2|_W \ge G(x, w_2) - G(x,w_1) \ge G(x, w_2) - 1 \ge \mu_2(r),
$$
and \eqref{ek} yields
\be{dc}
d_H(Z(w_1),Z(w_2)) \le D_K \le \frac{D_K L}{\mu_2(r)}|w_1 - w_2|_W.
\ee
In the case 
\\[2mm]
B. $d_H(Z(w_1),Z(w_2))<r$
\\[2mm]
we proceed as follows. For every $\ve > 0$ there exists 
$x_\ve \in Z(w_1)$ be such that $d_\ve :=\dist(x_\ve, Z(w_2)) \in (0,r)$ and 
$d_H(Z(w_1),Z(w_2)) - \ve \le d_\ve \le d_H(Z(w_1),Z(w_2))$.
By Proposition \ref{lp1}, there exists $x'_\ve \in \partial Z(w_2)$ such that 
$x_\ve = x'_\ve + d_\ve n'_\ve$, where
\begin{equation}\label{dc-2}
n'_\ve = \frac{\nabla_x G(x'_\ve,w_2)}{|\nabla_x G(x'_\ve,w_2)|}.
\end{equation}
We have $G(x'_\ve,w_2) = 1$, $G(x_\ve,w_2) > 1$, and 
by \eqref{dc-2}, Hypotheses \ref{hl}\,(i) and \ref{hl}\,(iii),
\begin{align*}
& G(x_\ve,w_2) - G(x'_\ve,w_2) = \int_0^1 \scal{\nabla_x G(x'_\ve+ td_\ve n'_\ve,w_2), d_\ve n'}\dd t\\
& = \scal{\nabla_x G(x'_\ve,w_2), d_\ve n'_\ve} + 
     \int_0^1 \scal{\nabla_x G(x'_\ve+ td_\ve n'_\ve,w_2)-\nabla_x G(x'_\ve,w_2), d_\ve n'_\ve}\dd t\\
&\ge cd_\ve - \frac{\lambda}{2} d_\ve^2 = \frac{c}{2}d_\ve\left(2 - \frac{d_\ve}{r}\right) \ge \frac{cd_\ve}{2}.
\end{align*}
On the other hand, we have $G(x_\ve,w_1) \le 1 = G(x'_\ve,w_2)$, hence, by Hypothesis \ref{hl}\,(iv),
\begin{align}\label{dd}
d_H(Z(w_1),Z(w_2)) - \ve & \le d_\ve \le \frac{2}{c}(G(x_\ve,w_2) - G(x'_\ve,w_2)) \notag \\
 & \le \frac{2}{c}(G(x_\ve,w_2) - G(x_\ve,w_1)) \le \frac{2L}{c}|w_1-w_2|_W,
\end{align}
and combining \eqref{dc} with \eqref{dd} and with the arbitrariness of $\ve$ we complete the proof.
\epf

We cite without proof the following result.

\begin{proposition}\label{ap1}
Let $\{Z(w); w \in W\}$ be a family of $r$-prox-regular sets and 
let \eqref{lipw} hold for every $K > 0$ and every $w_1,w_2 \in W$.
Then for every $u \in W^{1,1}(0,T; X)$, $w \in W^{1,1}(0,T; W)$, and every initial condition $x_0 \in Z(w(0))$ there exists a unique solution $\xi\in W^{1,1}(0,T; X)$ such that 
\begin{align}
 &\ \  x(t) := u(t) - \xi(t) \in Z(w(t)) \quad \text{for every $t \in [0,T]$}, \label{GP1}\\
 & \scal{x(t) - z, \dot\xi(t)} + \frac{|\dot\xi(t)|}{2r}|x(t) - z|^2 \ge 0 
   \quad \forall z \in Z(w(t)), \ \text{for a.e. $t \in (0,T)$}, \label{GP2}\\
 &\ \  x(0) = x_0. \label{GP3}
\end{align}
\end{proposition}

The statement was proved in \cite[Corollary 5.3]{kmr} under the assumption that the constant $C_K$ in \eqref{lipw} can be chosen independently of $K$. This is indeed not a real restriction, since the input values $w(t)$ belong to an a priori bounded set. We obtain global Lipschitz continuity under the hypotheses of Proposition \ref{ap1} by choosing
$K > \sup_{t\in [0,T]} |w(t)|_W$ in \eqref{lipw} and modifying the function $G$ for
$|w|_W \ge K$ for instance as $\tilde G(x,w) = G(x, f(|w|_W) w)$, where
$f:[0,\infty) \to [0,\infty)$ is a smooth function such that $f(s) = 1$
for $s\in[0,K]$ and $s f(s) \le K$ for $s> K$.

The above developments have shown that the assumptions of Proposition \ref{ap1} are fulfilled if Hypothesis \ref{hl} holds.
The existence and uniqueness of solutions to \eqref{a1} is therefore guaranteed for all $u \in W^{1,1}(0,T; X)$, $w \in W^{1,1}(0,T; W)$, and every initial condition $x_0 \in Z(w(0))$. We now prove the following identity which plays a substantial role in our arguments.

\begin{lemma}\label{al2}
Let Hypothesis \ref{hl} hold and let $u,w,\xi,x$ be as in Proposition \ref{ap1}. 
Let  $\nabla_w G: X \times W \to W'$ be continuous. Then
\be{a6}
\scal{\dot \xi(t),\dot x(t) + s(t)} = 0
\ee
for almost all $t \in (0,T)$ with the choice
\be{w2}
s(t) = \frac{\nabla_x G(x(t),w(t))}{\dist(x(t), \partial Z(w(t))+|\nabla_x G(x(t),w(t))|^2}\duw{\dot w(t), \nabla_w G(x(t),w(t))}.
\ee
\end{lemma}

\bpf{Proof}
We first check that the denominator in \eqref{w2} is bounded away from zero. Indeed, thanks to Hypothesis \ref{hl}(ii) we find $\delta_c > 0$ such that the implication
$$
|x_1 - x_2| < \delta_c \ \Longrightarrow \ |\nabla_x G(x_1,w(t))- \nabla_x G(x_2,w(t))| < \frac{c}{2} \quad \forall t\in [0,T]
$$
holds for all $x_1, x_2 \in Z(w(t))$. Then we have
$$
\dist(x(t), \partial Z(w(t))) < \delta_c  \, \Rightarrow \, \exists \hat x \in \partial Z(w(t)): |\hat x - x(t)| < \delta_c  \, \Rightarrow \, |\nabla_x G(x(t),w(t))| > \frac{c}{2}.
$$
With this choice of $\hat x$, we have
\be{low}
\dist(x(t), \partial Z(w(t)))+|\nabla_x G(x(t),w(t))|^2 \ge \min\left\{\delta_c, \frac{c^2}{4}\right\}
\ee
for all $t \in [0,T]$. For a.~e. $t \in (0,T)$ one of the following two cases occurs:
\begin{itemize}
\item[(1)] $\dot\xi(t) = 0$,
\item[(2)] $\dot\xi(t) \ne 0$.
\end{itemize}
Let $B = \{t \in (0,T): \dot\xi(t) \ne 0\}$.
For $t\in B$ we have $x(t) \in \partial Z(w(t))$, that is, $G(x(t),w(t)) = 1$. Hence, for a.~e. $t \in B$ we have
\be{acc1}
\scal{\dot x(t), \nabla_x G(x(t),w(t))} + \duw{\dot w(t), \nabla_w G(x(t),w(t))} = 0.
\ee
Moreover, for $t\in B$, the vector $\dot\xi(t)$ points in the direction of the unit outward normal vector $n(x(t),w(t))$ to $Z(w(t))$ at the point $x(t)$, that is, 
\be{w0}
\dot \xi(t) = \frac{|\dot \xi(t)|}{|\nabla_x G(x(t),w(t))|} \nabla_x G(x(t),w(t)).
\ee
From \eqref{acc1}--\eqref{w0} we obtain for $t \in B$ the identity
\be{acc2}
\scal{\dot \xi(t),\dot x(t)} + \scal{\dot \xi(t),\frac{\nabla_x G(x(t),w(t))}{|\nabla_x G(x(t),w(t))|^2}\duw{\dot w(t), \nabla_w G(x(t),w(t))}} = 0,
\ee
which is of the desired form
\be{acc3}
\scal{\dot \xi(t),\dot x(t) + \hat s(t)} = 0
\ee
with
\be{w3}
\hat s(t) = \frac{\nabla_x G(x(t),w(t))}{|\nabla_x G(x(t),w(t))|^2}\duw{\dot w(t), \nabla_w G(x(t),w(t))},
\ee
which holds for a.~e. $t \in B$ by the above argument. Since $\dist(x(t), \partial Z(w(t))) = 0$ for $t \in B$, we obtain \eqref{a6} directly from \eqref{acc3}. For $t \in (0,T)\setminus B$, identity \eqref{a6} is trivial since $\dot\xi(t) = 0$ a.\,e. on $t \in (0,T)\setminus B$.
\epf

Under Hypothesis \ref{hl},  the mapping $(x,w) \mapsto \dist(x,Z(w))$ is locally Lipschitz continuous. This can be easily proved as follows. Let $x,x' \in X$, $w,w' \in W$ be given. Put $d=\dist (x, Z(w))$, $d'=\dist (x', Z(w'))$, and assume for instance that $d \ge d'$. For an arbitrary $\ve > 0$ we find $z'\in Z(w')$ such that $|x'-z'|\le d'+\ve$, and $z \in Z(w)$ such that $|z-z'|\le d_H(Z(w),Z(w'))+\ve$. Then
\be{dis}
d \le |x-z| \le |x-x'|+ |x'-z'|+|z'-z| \le |x-x'|+  d' + d_H(Z(w),Z(w'))+2\ve
\ee
and the assertion follows from \eqref{lipw}. Moreover the following statement holds true.

\begin{lemma}\label{ldi}
Let Hypothesis \ref{hl} hold, and let $K>0$ be given. Then there exists $m_K>0$ such that for all $w,w' \in W$ satisfying the inequalities
\be{emstar}
\max\{|w|_W, |w'|_W\} \le K,\quad |w - w'|_W < m_K
\ee
and for all $x, x' \in X$, $x \in Z(w)$, $x' \in Z(w')$ we have
$$
|\dist(x,\partial Z(w)) - \dist(x',\partial Z(w'))| \le |x-x'| + d_H(Z(w), Z(w')).
$$
\end{lemma}

\bpf{Proof}
Put $d = \dist(x,\partial Z(w))$, $d' =\dist(x',\partial Z(w'))$, and assume $d \ge d'$. 
For every $\ve>0$ we find $z' \in \partial Z(w')$ such that $|x'-z'| \le d'+\ve$, and $z \in \partial Z(w)$ such that $|z'-z| \le d_H(\partial Z(w), \partial Z(w'))+\ve$. We argue as in \eqref{dis} and obtain
$$
d \le |x-z| \le |x-x'| + |x' - z'| + |z' - z| \le |x-x'| + d' + d_H(\partial Z(w), \partial Z(w'))+2\ve.
$$
The proof will be complete if we prove that for a suitable value of $m_K$ and for $w,w'$ satisfying \eqref{emstar} we have
\be{hd}
d_H(\partial Z(w),\partial Z(w')) \le d_H(Z(w), Z(w')).
\ee
We claim that the right choice of $m_K$ is
$$
m_K = \frac{d^*}{C_K}
$$
with $C_K$ from Lemma \ref{ldh} and any $d^* < r$ with $r$ as in Proposition \ref{lp1}.

Indeed, from Lemma \ref{ldh} it follows that $\rho :=$ 
$d_H(Z(w), Z(w')) \le d^*$. Consider any $z \in \partial Z(w)$ and assume that $z \notin \partial Z(w')$. We distinguish two cases: $z \in \Int Z(w')$ and $z \notin Z(w')$. For $z \in \Int Z(w')$ and $t\ge 0$ we put
$$
z(t) = z+ t \frac{\nabla_x G(z,w)}{|\nabla_x G(z,w)|}.
$$
Then for $t<r$ we have $\dist(z(t), Z(w)) = |z(t) - z| = t$. Since $d_H(Z(w), Z(w')) =\rho \le d^*< r$, there exists necessarily $t \le \rho$ such that $z(t) \in \partial Z(w')$, and we conclude that 
$\dist(z,\partial Z(w')) \le \rho$. In the case $z \notin Z(w')$ we use Proposition \ref{lp1} and find 
$z' \in \partial Z(w')$ such that $\dist(z, Z(w')) = \dist(z,\partial Z(w')) = |z-z'| \le \rho$ and \eqref{hd} follows.
\epf

It is easy to see that a counterpart of inequality \eqref{hd} does not hold for general sets. It suffices to consider $R_1 > R_2 > 0$ and $Z_1 = \overline{B_{R_1}(0)}$, $Z_2 = Z_1 \setminus B_{R_2}(0)$,  where for 
$x \in X$ and $R > 0$ we denote by $B_R(x)$ the open ball $\{y\in X: |x-y| < R\}$. Then 
$d_H(Z_1, Z_2) = R_2$, $d_H(\partial Z_1, \partial Z_2) = R_1-R_2$, so that \eqref{hd} is violated for 
$R_2 < R_1/2$.

The solution mapping of \eqref{GP1}--\eqref{GP3} is continuous in the following sense.

\begin{theorem}\label{t3}
Let Hypothesis \ref{hl} hold and let $\nabla_w G: X \times W \to W'$ be continuous.
 Let $u\in W^{1,1}(0,T; X)$ and $w \in W^{1,1}(0,T; W)$ be given, and let $\{u_n; n\in \nat\} \subset W^{1,1}(0,T; X)$ and $\{w_n; n\in \nat\} \subset W^{1,1}(0,T; W)$ be sequences such that $u_n(0) \to u(0)$, $w_n(0) \to w(0)$ as $n \to \infty$, and
\be{a8}
\lim_{n\to \infty}\int_0^T \left(|\dot u_n(t) - \dot u(t)| + |\dot w_n(t) - \dot w(t)|_W\right)\dd t = 0
\ee
as $n \to \infty$. Let $\xi_n, \xi \in W^{1,1}(0,T; X)$ be the solutions to \eqref{GP1}--\eqref{GP3}
corresponding to the inputs $u_n, w_n, u,w$, respectively, with initial conditions 
$x_n^0 \in Z(w_n(0)), x_0 \in Z(w(0))$ such that $|x_n^0 - x_0| \to 0$ as $n \to \infty$. Then
\be{a9}
\lim_{n\to \infty}\int_0^T |\dot \xi_n(t) - \dot \xi(t)|\dd t = 0.
\ee
\end{theorem}

The proof of Theorem \ref{t3} relies on the following general property of functions in $L^1(0,T;X)$ proved in \cite{vhm}.

\begin{lemma}\label{al3}
Let $\{v_n; n\in \nat \cup \{0\}\} \subset
L^1(0,T;X)$, $\{g_n; n\in \nat\cup \{0\}\} \subset L^1(0,T; \real)$ be given sequences such that
\begin{itemize}
\item[{\rm (i)}]\, $\lim_{n\to \infty} \int^T_0\scal{v_n(t), \varphi(t)} \dd t = \int^T_0\scal{v(t), \varphi(t)}\dd t
 \quad \forall \varphi\in C([a,b];X)$, 
\item[{\rm (ii)}]\, $\lim_{n\to \infty} \int^T_0|g_n(t) - g_0(t)| \dd t = 0$,
\item[{\rm (iii)}]\, $|v_n(t)| \le g_n(t)$ a.\,e.\, $\forall n\in \nat$,
\item[{\rm (iv)}]\, $|v_0(t)| = g_0(t)$ a.\,e.
\end{itemize}
Then $\lim_{n\to \infty} \int^T_0|v_n(t) - v_0(t)| \dd t = 0$.
\end{lemma}

Notice that Lemma \ref{al3} does not follow from the Lebesgue Dominated Convergence Theorem,
since {\it we do not assume the pointwise convergence\/}. The proof is elementary and we repeat
it here for the reader's convenience.

\bpf{Proof of Lemma \ref{al3}}
We first prove that property (i) holds for every $\varphi \in L^\infty (0,T; X)$. For a fixed
$\varphi\in L^\infty(0,T;X)$ and $\delta > 0$ we use Lusin's Theorem
to find a function $\psi \in C([0,T];X)$\, and a set\, $M_\delta \subset
[0,T]$ such that $\meas (M_\delta) < \delta$ and $\psi(t) = \varphi(t)$
for all $t\in [0,T] \setminus M_\delta, \|\psi\| \le
\|\varphi\|$. We then have 
$$
\begin{aligned}
\left|\int_0^T\scal{v_n(t) - v_0(t), \varphi(t)} \dd t\right| &\le \left|
\int_0^T \scal{v_n(t) - v_0(t), \psi(t)}\dd t \right| \\
&+ 2\|\varphi\| \left(\int_0^T |g_n(t)-g_0(t)| \dd t +
2\int_{M_\delta}g_0(t) \dd t\right).
\end{aligned}
$$
Since $\delta$ can be chosen arbitrarily small and $g_0\in L^1(0,T)$, the integral of $g_0$ over $M_\delta$ can be made arbitrarily small and we obtain 
\be{a10}
\lim_{n\to \infty} \int_0^T \scal{v_n (t), \varphi(t)} \dd t =
\int_0^T\scal{v_0(t), \varphi(t)} \dd t \quad \forall\varphi \in
L^\infty(0,T;X).
\ee
Let us note that the transition from (i) to \eqref{a10} is related to
the Dunford-Pettis Theorem, see \cite{edw}. To prove Lemma \ref{al3} we
put for $t\in [0,T]$
$$
\varphi(t) := \left\{
\begin{array}{ll} 
0&\text{if}\quad v_0(t) = 0,\\
\frac{v_0(t)}{g_0(t)} \quad &\text{if} \quad v_0(t) \ne 0.
\end{array}
\right.
$$
Then $\varphi\in L^\infty(0,T;X)$ and the inequality 
$$
\begin{aligned}
|v_n(t) -v_0(t)|^2 &\le g^2_n(t) - 2\scal{v_n(t), v_0(t)} + g^2_0(t)\\
&= |g_n(t) - g_0(t)|^2 + 2g_0(t)\big(g_n(t) - g_0(t) +
\scal{v_0(t), \varphi(t)} - \scal{v_n(t),\varphi(t)}\big) 
\end{aligned}
$$
holds for a.\,e. $t\in [0,T]$. By H\"older's inequality we have 
$$
\begin{aligned}
&\int_0^T|v_n(t) - v_0(t)| \dd t \le \int_0^T|g_n(t) -
g_0(t)|\dd t\\
&\qquad +\left(\int_0^T 2g_0(t)\dd t\right)^{1/2}
\left(\int_0^T\big(g_n(t)- g_0(t) + \scal{v_0(t), \varphi(t)} -
\scal{v_n(t), \varphi(t)} \big) \dd t \right)^{1/2},
\end{aligned}
$$
and the assertion follows from \eqref{a10}.
\epf

We are now ready to prove one of our main results, namely Theorem \ref{t3}.

\bpf{Proof of Theorem \ref{t3}}
By Lemma \ref{al2} we check that $s_n$ given by the formula
\be{w4}
s_n(t) = \frac{\nabla_x G(x_n(t),w_n(t))}{\dist(x_n(t), \partial Z(w_n(t)))+|\nabla_x G(x_n(t),w_n(t))|^2}\duw{\dot w_n(t), \nabla_w G(x_n(t),w_n(t))}
\ee
satisfy a.\,e. the identity 
\be{w5}
\scal{\dot \xi_n(t),\dot x_n(t) + s_n(t)} = 0.
\ee
Theorem 4.4 of \cite{kmr} states that $x_n \to x$ uniformly in $C([0,T];X)$. Using Lemma~\ref{ldi} and formulas \eqref{low}, \eqref{a8} we conclude that $s_n$ converge strongly to $s$ in $L^1(0,T;X)$.
Put $y_n = \dot u_n +s_n - 2\dot\xi_n = \dot x_n + s_n - \dot \xi_n$, $y = \dot u +s - 2\dot\xi = \dot x + s - \dot \xi$. For a.\,e. $t\in (0,T)$ we have by \eqref{a6}, \eqref{w5} that
\be{a11}
|y_n(t)|^2 = |\dot x_n(t) + s_n(t) - \dot \xi_n(t)|^2 = |\dot x_n(t) + s_n(t) + \dot \xi_n(t)|^2 = |\dot u_n(t) + s_n(t)|^2,
\ee
and similarly $|y(t)|^2 = |\dot u(t) + s(t)|^2$. Put $v_n(t) = y_n(t)$, $v_0(t) = y(t)$, $g_n(t) = |\dot u_n(t) + s_n(t)|$, $g_0(t) = |\dot u(t) + s(t)|$. 
We see that hypotheses (ii)--(iv) of Lemma \ref{al3} are satisfied. The assertion of Theorem \ref{t3} will follow from Lemma \ref{al3} provided we check that
\be{a12}
\lim_{n\to \infty} \int_0^T \scal{y_n (t) - y(t), \varphi(t)} \dd t = 0 \quad \forall\varphi \in C([0,T];X).
\ee
This will certainly be true if we prove that
\be{w6}
\lim_{n\to \infty} \int_0^T \scal{\dot \xi_n (t) - \dot \xi(t), \varphi(t)} \dd t = 0 \quad \forall\varphi \in C([0,T];X).
\ee
Let $\varphi \in C([0,T];X)$ be given. For an arbitrary $\ve >0$ we find $\psi \in C^1([0,T];X)$ such that $\|\psi - \varphi\| < \ve$. There exists a constant $C>0$ independent of $n$ and $\ve$ such that
$$
\left|\int_0^T \scal{\dot \xi_n (t) - \dot \xi(t), \varphi(t)- \psi(t)} \dd t\right| \le C\ve,
$$
hence,
\be{a13a}
\left|\int_0^T \scal{\dot \xi_n (t) - \dot \xi(t), \varphi(t)} \dd t\right| \le \left|\int_0^T \scal{\dot \xi_n (t) - \dot \xi(t), \psi(t)} \dd t\right| + C\ve,
\ee
where we can integrate by parts and obtain
\be{a13}
\begin{aligned}
\int_0^T \scal{\dot \xi_n (t) - \dot \xi(t), \psi(t)} \dd t &= \scal{\xi_n (T) - \xi(T), \psi(T)} - \scal{\xi_n (0) - \xi(0), \psi(0)}\\
&\quad - \int_0^T \scal{\xi_n (t) - \xi(t), \dot\psi(t)} \dd t.
\end{aligned}
\ee
By \cite[Theorem 4.4]{kmr}, the right-hand side of \eqref{a13} converges to $0$ as $n \to \infty$. Since $\ve$ in \eqref{a13a} can be chosen arbitrarily small, we obtain \eqref{a12} from \eqref{a13a} and \eqref{a13}.
Using Lemma \ref{al3} we conclude that
$$
\lim_{n\to \infty}\int_0^T |y_n(t) - y(t)|\dd t = 0
$$
and the assertion of Theorem \ref{t3} easily follows.
\epf


\section{Local Lipschitz continuity}\label{lip}

We have proved in the previous section that the solution mapping $(u,w) \mapsto \xi$ of Problem 
\eqref{GP1}--\eqref{GP3} is strongly continuous with respect to the $W^{1,1}$-norm provided Hypothesis \ref{hl} holds and $\nabla_{w}G$ is a continuous function. 
Here we show that if $\nabla_{x}G$, $\nabla_{w}G$ are Lipschitz continuous, then the solution mapping of Problem \eqref{GP1}--\eqref{GP3} is locally Lipschitz continuous with respect to the $W^{1,1}$-norm. 
Here are the precise assumptions.

\begin{hyp}\label{imh1}
Let Hypothesis \ref{hl} hold. Assume that the partial derivatives $\nabla_{x}G(x,w)$ $\in$ $X$, 
$\nabla_w\,G(x,w) \in W'$ exist for every $(x,w) \in X\times W$ and
there exist positive constants $K_0, K_1, C_0, C_1$ such that

\begin{itemize}
\item[{\rm (i)}]\
$|\nabla_x G(x,w)| \le K_0, \ |\nabla_w G(x,w)|_{W'} \le K_1 \qquad \forall (x,w) \in X\times W$,

\item[{\rm (ii)}]\
for every $(x,w), (x',w') \in X\times W$ we have
\bearr\label{2e6}
|\nabla_{x}G(x,w)-\nabla_{x}G(x',w')| &\le& C_0 (|x-x'|+|w-w'|_W)\,,
\\[1mm] \label{2e7}
|\nabla_{w}G(x,w)-\nabla_{w}G(x',w')|_{W'} &\le& C_1\, (|x-x'|+|w-w'|_W).
\eearr
\end{itemize}
\end{hyp}

In the following two lemmas we derive some useful formulas.

\begin{lemma}\label{2l2}
Let Hypothesis \ref{imh1}\,(i) hold, and let $(u,w) \in \domab$, $x^0 \in Z({w}(0))$,  
and $\xi \in \AC$ satisfy \eqref{GP1}--\eqref{GP3} with $x(t) = u(t) - \xi(t)$. 
For $t\in (0,T)$ set
\bears
A[u,w](t) &=& \scal{\doxi(t), \nabla_{x}G(x(t),w(t))}\, ,\\[1mm]
B[u,w](t) &=& \scal{\dou(t), \nabla_{x}G(x(t),w(t))} +
\duw{\dot w(t),\nabla_{w}G(x(t),w(t))}\,,
\eears
Then for a.\,e. $t \in (0,T)$ we have either
\begin{itemize}
\item[{\rm (i)}] $\doxi(t) = 0$, $\frac{\dd}{\dd t} G(x(t),w(t)) = B[u,w](t)$,
\item[{}] \qquad or
\item[{\rm (ii)}] $\doxi(t) \ne 0$, $x(t) \in \partial Z(w(t))$,
$A[u,w](t) = B[u,w](t) > 0$, $\max_{\tau\in [0,T]} G(x(\tau),w(\tau)) = G(x(t),w(t)) = 1$, $\frac{\dd}{\dd t} G(x(t),w(t)) = 0$, and
\be{2e7a}
\doxi(t) = \frac{A[u,w](t)}{|\nabla_{x}G(x(t),w(t))|^2}\,\nabla_{x}G(x(t),w(t))\,.
\ee
\end{itemize}
Moreover, for a.\,e. $t \in (0,T)$ we have that
\bearr\label{2e8b}
|B[u,w](t)| &\le& |\nabla_x G(x(t),w(t))|\, |\dou(t)| + K_1 |\dow(t)|_W ,
\\[1mm] \label{2e8a}
|\doxi(t)| &\le& |\dou(t)| +  \frac{K_1}{c} |\dow(t)|_W,
\eearr
with $c$ from Hypothesis \ref{hl}\,(i).
\end{lemma}

\bpf{Proof}\
Let $L \subset (0,T)$ be the set of Lebesgue points of all functions
$\dou$, $\dot w$, $\doxi$. Then $L$ has full measure
in $[0,T]$, and for $t\in L$ we have
\be{2e8}
\frac {\dd}{\dd t}G(x(t),w(t))  = \scal{\dox(t), \nabla_{x}G(x(t),w(t))} +
\duw{\dot w(t),\nabla_{w}G(x(t),w(t))}\, .
\ee
If $\doxi(t) = 0$, then $\dox(t) = \dou(t)$, and (i) follows from
(\ref{2e8}). If $\doxi(t) \ne 0$, then $x(t) \in \partial Z({w}(t))$,
hence $G(x(t),w(t)) = 1 = \max_{\tau\in [0,T]} G(x(\tau),w(\tau))$ and $\frac {\dd}{\dd t}G(x(t),w(t))=0$, so that (\ref{2e7a})
follows from \eqref{w0}. Furthermore,
(\ref{2e8}) yields $\scal{\dox(t), \nabla_{x}G(x(t),w(t))} = -\duw{\dot w(t),\nabla_{w}G(x(t),w(t))}$, hence
\bears
\scal{\doxi(t), \nabla_{x}G(x(t),w(t))}
&=& \scal{\dou(t), \nabla_{x}G(x(t),w(t))} - \scal{\dox(t), \nabla_{x}G(x(t),w(t))}\\[1mm]
&=& \scal{\dou(t), \nabla_{x}G(x(t),w(t))} + \duw{\dot w(t),\nabla_{w}G(x(t),w(t))}\, ,
\eears
We are left to prove \eqref{2e8b}--\eqref{2e8a}.
Formula \eqref{2e8b} follows from Hypothesis \ref{imh1}\,(i). Formula
(\ref{2e8a}) is trivial if $\doxi(t) = 0$; otherwise we
have $|\doxi(t)|$ $=$ $A[u,w](t)/|\nabla_{x}G(x(t),w(t))|$ $=$ $B[u,w](t)/|\nabla_{x}G(x(t),w(t))|$, 
$x(t) \in \partial Z(w(t))$, and (\ref{2e8a}) follows.
\epf

\begin{lemma}\label{2l3}
Let Hypothesis \ref{imh1}\,(i) hold, let $({u_i},{w_i}) \in \domab$ and $x_i^0 \in Z({w_i}(0))$ be given for $i=1,2$,
let $\xi_i \in \AC$ be the respective solutions to 
\eqref{GP1}--\eqref{GP3} with
$x_i = {u_i}-\xi_i$ for $i=1,2$.
Then for a.\,e. $t\in (0,T)$ we have
\bearr\label{2e9}
|A[u_1,w_1](t) {-} A[u_2,w_2](t)| &+&
\frac {\dd}{\dd t} |G(x_1(t),w_1(t))-G(x_2(t),w_2(t))|\\ \nonumber
&\le& |B[u_1,w_1](t) - B[u_2,w_2](t)|\, ,\\[1mm] \label{2e10}
|\doxi_1(t) - \doxx(t)| &\le& \frac1c\,|A[u_1,w_1](t) - A[u_2,w_2](t)|
\\[1mm] \nonumber
&&\hspace{-24mm} +\ \frac1c\,\left(|\dou_1(t)| +  \frac{K_1}c |\dot w_1(t)|_W\right)
|\nabla_x G(x_1(t),w_1(t))-\nabla_x G(x_2(t),w_2(t))|\, .
\eearr
where $A$ and $B$ are defined as in Lemma \ref{2l2}.
\end{lemma}

\bpf{Proof}\
The assertion follows directly from Lemma \ref{2l2} if $\doxi_1(t) = \doxx(t) = 0$. Assume now
\smallskip

$\bullet$\ $\doxi_1(t)\ne 0$, $\doxx(t) \ne 0$.

Then (\ref{2e9}) is again an immediate consequence of Lemma \ref{2l2}.
To prove  (\ref{2e10}), we use (\ref{2e7a}) and the elementary
vector identity
\[
\left|\frac{z}{|z|^2} - \frac{z'}{|z'|^2}\right| \ = \ \frac 1{|z| |z'|}\,
|z - z'| \quad \for \ z, z' \in X \setminus\{0\}\, ,
\]
to obtain
\bears
|\doxi_1(t) - \doxx(t)| &\le& |A[u_1,w_1](t)|\,
\left|\frac{\nabla_x G(x_1(t),w_1(t))}{|\nabla_x G(x_1(t),w_1(t))|^2} -
\frac{\nabla_x G(x_2(t),w_2(t))}{|\nabla_x G(x_2(t),w_2(t))|^2}\right|\\[1mm]
&& + \ \frac{1}{|\nabla_x G(x_2(t),w_2(t))|}\,|A[u_1,w_1](t) - A[u_2,w_2](t)|\\[1mm]
&& \hspace{-21mm}=\
\frac{|B[u_1,w_1](t)|}{|\nabla_x G(x_1(t),w_1(t))|\,|\nabla_x G(x_2(t),w_2(t))|} \left|\nabla_x G(x_1(t),w_1(t)){-}\nabla_x G(x_2(t),w_2(t))\right|\\[1mm]
&& + \ \frac{1}{|\nabla_x G(x_2(t),w_2(t))|}\,|A[u_1,w_1](t) - A[u_2,w_2](t)| .
\eears
By Hypothesis \ref{hl}\,(i) we have $|\nabla_x G(x_i(t),w_i(t))| \ge c$ for $i=1,2$, and
combining the above inequalities with (\ref{2e8b}) we obtain the assertion.

Let us consider now the case
\smallskip

$\bullet$\ $\doxi_1(t)  \ne 0$, $\doxx(t) = 0$.

Then $|A[u_1,w_1](t) - A[u_2,w_2](t)| = A[u_1,w_1](t)$,
$G(x_1(t),w_1(t))-G(x_2(t),w_2(t)) = 1 - G(x_2(t),w_2(t)) \ge 0$, hence
\bears
&&\hspace{-20mm}|A[u_1,w_1](t) - A[u_2,w_2](t)| + \frac \dd{\dd t}
|G(x_1(t),w_1(t))-G(x_2(t),w_2(t))|\\[1mm]
&=& A[u_1,w_1](t) - \frac \dd{\dd t} G(x_2(t),w_2(t))\\[1mm]
&=& B[u_1,w_1](t) - B[u_2,w_2](t)\, ,
\eears
hence (\ref{2e9}) is fulfilled. We further have similarly as above
that
\[
|\doxi_1(t) - \doxx(t)|\, =\, |\doxi_1(t)| \le \frac1c\, A[{v_1}, {u_1}](t)
= \frac1c\, |A[{u_1}, {w_1}](t) - A[{u_2}, {w_2}](t)|\, ,
\]
hence (\ref{2e10}) holds.
The remaining case
\smallskip

$\bullet$\ $\doxi_1(t) = 0$, $\doxx(t) \ne 0$

\noindent is analogous, and Lemma \ref{2l3} is proved.
\epf

We are now ready to prove the following main result.

\begin{theorem}\label{2p4}
Let Hypothesis \ref{imh1} hold, let $({u_i},{w_i})\in \domab$
and $x_i^0 \in Z({w_i}(0))$ be given for $i=1,2$, let $\xi_i \in \AC$ be the respective solutions to
\eqref{GP1}--\eqref{GP3} with $x_i = {u_i}-\xi_i$ for $i=1,2$. Then for a.\,e. $t\in (0,T)$ we have
\begin{align}\nonumber
|\doxi_1(t) - \doxx(t)| &+ \frac1c\,
\frac \dd{\dd t} |G(x_1(t),w_1(t))-G(x_2(t),w_2(t))|\\[1mm] \nonumber 
&\quad \le \frac{1}{c} \big(K_0|\dou_1(t) - \dou_2(t)| + K_1 |\dot w_1(t) - \dot w_2(t)|_W\big)
\\[1mm]\label{2e11}
&\hspace{-22mm} + \frac1c \left(2C_0 |\dou_1(t)| +
\left(C_1 +\frac{C_0 K_1}{c}\right)\,|\dot w_1(t)|_W \right)\big(|{w_1}(t) {-} {w_2}(t)|_W + |{x_1}(t) {-} {x_2}(t)|\big) .
\end{align}
\end{theorem}

\bpf{Proof}\
By Lemma \ref{2l3}, we have
\begin{align}
  & |\doxi_1(t) - \doxx(t)| + 
      \frac1c\,\frac \dd{\dd t} |G(x_1(t),w_1(t))-G(x_2(t),w_2(t))|
       \le \frac1c|B[u_1,w_1](t) - B[u_2,w_2](t)| \notag \\
  & \quad +\ \frac1c\,\left(|\dou_1(t)| +  \frac{K_1}c |\dot w_1(t)|_W\right)
|\nabla_x G(x_1(t),w_1(t))-\nabla_x G(x_2(t),w_2(t))|, \notag
\end{align}
where $B[u,w]$ is defined as in Lemma \ref{2l2}. Hence 
(\ref{2e11}) follows since
from Hypothesis \ref{imh1} and the triangle inequality applied to $B[u,w]$ 
we infer that
\begin{align*}
&|B[u_1,w_1](t) - B[u_2,w_2](t)| \le   K_0 |\dot u_1(t) - \dot u_2(t)|  + K_1|\dot w_1(t) - \dot w_2(t)|_W\\
   & \qquad +
	\big(C_0|\dot u_1(t)|+ C_1|\dot w_1(t)|\big)\big(|x_1(t) - x_2(t)| + |w_1(t) - w_2(t)|_W\big).
\end{align*}
\epf

\begin{corollary}\label{2c5}
For every $R>0$ there exists a constant $C(R)>0$ such that for every $({u_i},{w_i})\in \domab$ and every $x_i^0 \in Z({w_i}(0))$ for $i=1,2$ such that $\max\{\int_0^T|\dot u_1(t)|\dd t, \int_0^T|\dot w_1(t)|_W\dd t\} \le R$, the respective solutions $\xi_i \in \AC$ to \eqref{GP1}--\eqref{GP3} 
satisfy the inequality
\begin{align}\nonumber
\int_0^T|\dot \xi_1(t) - \dot \xi_2(t)|\dd t 
&\le C(R) \bigg(\int_0^T\left(|\dot u_1(t) - \dot u_2(t)| + |\dot w_1(t) - \dot w_2(t)|_W\right)\dd t \\ \label{e3}
& \quad +|w_1(0) - w_2(0)|_W  + |x_1^0 - x_2^0|\bigg).
\end{align}
\end{corollary}

\bpf{Proof}
In the situation of Theorem \ref{2p4} put $K_2 = \max\{K_0, K_1\}/c$, and
\begin{align*}
\Delta\xi(t) &= \int_0^t|\doxi_1(\tau) - \doxx(\tau)| \dd\tau + \frac1c\,
|G(x_1(t),w_1(t))-G(x_2(t),w_2(t))|,\\[1mm]
D(t) &= |\dou_1(t) - \dou_2(t)| + |\dot w_1(t) - \dot w_2(t)|_W,\\[1mm]
m(t) &= \frac1c \left(2C_0 |\dou_1(t)| +
\left(C_1 +\frac{C_0 K_1}{c}\right)\,|\dot w_1(t)|_W \right).
\end{align*}
Then from \eqref{2e11} it follows that
\begin{align}\nonumber
\frac \dd{\dd t}\Delta\xi(t) 
&\le K_2 D(t) + m(t) \big(|{w_1}(t) - {w_2}(t)|_W + 
   |{x_1}(t) - {x_2}(t)|\big)\\[1mm] \nonumber
&\le K_2 D(t) + m(t) \bigg(\int_0^t \big(|\dot w_1(\tau) - \dot w_2(\tau)|_W + 
   |\dot x_1(\tau) - \dot x_2(\tau)|\big)\dd\tau\\[1mm] \nonumber
&\qquad  + |w_1(0) - w_2(0)|_W  + |x_1^0 - x_2^0|\bigg)\\[1mm] \label{e4}
&\le K_2 D(t) + m(t) \bigg(\Delta\xi(t)
+ \int_0^t D(\tau)\dd\tau+ |w_1(0) {-} w_2(0)|_W  + |x_1^0 {-} x_2^0|\bigg).
\end{align}
We now use Gronwall's argument and put $M(t) = \int_0^t m(\tau)\dd\tau$. It follows from \eqref{e4} that
$$
\frac \dd{\dd t}\left(\expe^{-M(t)}\Delta\xi(t)\right)
\le \expe^{-M(t)}\bigg(K_2 D(t) + m(t) \bigg(\int_0^t D(\tau)\dd\tau+ |w_1(0) {-} w_2(0)|_W  + |x_1^0 {-} x_2^0|\bigg)\bigg).
$$
Integrating from $0$ to $T$ we obtain the assertion.
\epf

\begin{remark}
The local Lipschitz continuity of the input-output mapping cannot be expected if $\nabla_x G$ is not Lipschitz even if $G$ is convex. A counterexample is constructed in \cite[Theorem 2.2]{Kre01}. On the other hand, the global $W^{1,1}$-Lipschitz continuity of the sweeping process holds if $Z(w)$ is a convex polyhedron. This is shown for instance in \cite{DesTur99} for $Z$ independent of $w$, and it is generalized in \cite{KreVla01} to the case of non-orthogonal projections to the convex polyhedron, i.e. when the time derivative of $\xi$ lies in a prescribed cone of admissible directions.
\end{remark}


\section{Implicit sweeping processes}\label{imp}

In this section we consider the state dependent problem corresponding to \eqref{GP1}--\eqref{GP3},
where it is assumed that there exists a function $g: [0,T] \times X \times X \to W$ such that
\be{im1}
w(t) = g(t,u(t),\xi(t)).
\ee
More specifically,, given $u \in W^{1,1}(0,T; X)$ and $x_0 \in X$ such that
$x_0 \in Z(g(0,u(0),u(0)-x_0))$, one has to find $\xi \in W^{1,1}(0,T; X)$ such that
\begin{align}
 &\ \  x(t) := u(t) - \xi(t) \in Z(g(t,u(t),\xi(t))) \quad \text{for every $t \in [0,T]$}, \label{SGP1}\\
 & \scal{x(t) - z, \dot\xi(t)} + \frac{|\dot\xi(t)|}{2r}|x(t) - z|^2 \ge 0 
   \quad \forall z \in Z(g(t,u(t),\xi(t))) \ \text{ for a.\,e. } t \in (0,T), \label{SGP2}\\
 &\ \  x(0) = x_0. \label{SGP3}
\end{align}
Using the Banach contraction principle, we prove that \eqref{SGP1}--\eqref{SGP3}
 is uniquely solvable in $W^{1,1}(0,T; X)$ under the following assumptions on the function $g$.

\begin{hyp}\label{gh1}
A continuous function $g : [0,T] \times X \times X \to W$ is given such that its partial derivatives 
${\partial_t g, \partial_u g, \partial_\xi g}$ {exist and} satisfy the inequalities
\begin{align}\label{3e1}
|\partial_\xi g(t,u,\xi)|_{\L(X,W)} &\le \gamma\, ,
\\[1mm] \label{3e1a}
|\partial_u g(t,u,\xi)|_{\L(X,W)} &\le {\omega}\, ,
\\[1mm] \label{3e1au}
|\partial_t g(t,u,\xi)|_W &\le a(t)\, ,\\[1mm] \label{3e1b}
|\partial_\xi g(t,u,\xi) - \partial_\xi g(t,v,\eta)|_{\L(X,W)}
&\le C_\xi\,(|u - v| + |\xi - \eta|)
\, ,\\[1mm] \label{3e1c}
|\partial_u g(t,u,\xi) - \partial_u g(t,v,\eta)|_{\L(X,W)}
&\le C_u\,(|u - v| + |\xi - \eta|)
\, ,\\[1mm] \label{3e1cu}
|\partial_t g(t,u,\xi) - \partial_t g(t,v,\eta)|_W
&\le b(t)\,(|u - v| + |\xi - \eta|)
\end{align}
for every $u,v,\xi, \eta \in X$ and a.\,e. $t\in (0,T)$
with given functions $a, b\in L^1(0,T)$ and given
constants $\gamma, {\omega}, C_\xi,C_u >0$ such that
\be{3e2}
\delta := \frac{K_1\gamma}{c} < 1,
\ee
where $c$, $K_1$ are as in Hypothesis \ref{hl}\,(i) and Hypothesis \ref{imh1}.
\end{hyp}

Let us start our analysis with two auxiliary results.

\begin{lemma}\label{3l1}
Let Hypotheses \ref{imh1} and \ref{gh1} hold and let $\xi \in \AC$ 
satisfy \eqref{SGP1}--\eqref{SGP3}
with some $u \in \AC$ and some $x_0 \in X$ such that $x_0 \in Z(g(0,u(0),u(0)-x_0))$.
Then we have
\be{3e0}
|\doxi(t)| \le \frac 1{1-\delta}\, \left(\left(1 + \frac{\omega K_1}{c}\right) |\dou(t)| +  \frac{K_1}{c}\,a(t)\right)\quad \mbox{for a.\,e. $t \in (0,T)$.}
\ee
\end{lemma}

\bpf{Proof}
We set $w(t) = g(t,u(t),\xi(t))$ for $t \in [0,T]$. From Hypothesis \ref{gh1} follows that 
$w \in W^{1,1}(0,T; W)$, thus $u, w, \xi$ and $x_0$ satisfy \eqref{GP1}--\eqref{GP2} and Lemma \ref{2l2} applies. In particular
inequality (\ref{3e0}) is an easy consequence of \eqref{2e8a}. 
Indeed, using
(\ref{3e1}), (\ref{3e1a}) we obtain 
$|\dow(t)|_W \le a(t) +{\omega} |\dou(t)| + \gamma |\doxi(t)|$
so that (\ref{3e0}) follows from \eqref{3e2}.
\epf

Motivated by \eqref{3e0}, we define for any $u \in \AC$ the set
\be{ou}
\Omega(u):= \left\{\eta \in \AC\,:
\barr{rcl}
 |\dota(t)| &\le&
\frac 1{1-\delta}\, \left(\left(1 + \frac{\omega K_1}{c}\right) |\dou(t)| {+} \frac{K_1}{c}\,a(t)\right)\
\mbox{a.\,e.}\\[1mm]
\eta(0) &=& u(0)-x_0
\earr
\right\}
\ee
and prove the following statement.

\begin{lemma}\label{3l2}
For all $u \in \AC$ and $\eta\in \Omega(u)$, the solution $\xi \in \AC$ of \eqref{GP1}--\eqref{GP3} with $w(t) = g(t,u(t),\eta(t))$ belongs to $\Omega(u)$. Moreover, there exist constants $m_0 > 0$, $m_1 > 0$ such that for every $u_1, u_2 \in \AC$ and every $\eta_i \in \Omega(u_i)$, $i=1,2$, the solutions $\xi_i \in \AC$ of \eqref{GP1}--\eqref{GP3} with $w_i(t) = g(t,u_i(t),\eta_i(t))$ and $x_0^i \in X$ with 
$x_0^i \in Z(g(0,{u_i(0),} u_i(0)-x_0^i))$ $i=1,2$, satisfy for a.\,e. $t \in (0,T)$ the inequality
\begin{align}\nonumber
& |\doxi_1(t) - \doxx(t)| 
   + \frac1c\, \frac \dd{\dd t} |G(x_1(t),w_1(t))-G(x_2(t),w_2(t))| \notag \\
& \le m_1|\dou_1(t) - \dou_2(t)| + \delta |\dot \eta_1(t) - \dot \eta_2(t)| \notag
\\[1mm]\label{eta} 
& \quad + m_0\left(a(t) + b(t) + |\dou_1(t)|\right)\big(|{u_1}(t) - {u_2}(t)| +
|{\xi_1}(t) - {\xi_2}(t)| + |\eta_1(t) - \eta_2(t)|\big).
\end{align}
\end{lemma}

\bpf{Proof}
$u \in \AC$ and $\eta\in \Omega(u)$ be given. By \eqref{2e8a}, \eqref{im1}--\eqref{3e1au} and \eqref{3e2}, we have for a.\,e. $t \in (0,T)$ that
$$
|\dot\xi(t)| \le \left(1 + \frac{\omega K_1}{c}\right)|\dot u(t)| + \frac{K_1}{c}\left(a(t) + \gamma |\dot\eta(t)|\right)\le \frac 1{1{-}\delta}\, \left(\left(1 + \frac{\omega K_1}{c}\right) |\dou(t)| +  \frac{K_1}{c}\,a(t)\right),
$$
hence $\xi \in \Omega(u)$. To prove \eqref{eta}, we notice that the inequalities
\begin{align*}
|\dot w_1(t) - \dot w_2(t)|_W &\le \omega |\dot u_1(t) - \dot u_2(t)| + \gamma |\dot \eta_1(t) - \dot \eta_2(t)|\\[1mm]
&\quad + \left(b(t)+ C_u |\dot u_1(t)| + C_\xi |\dot \eta_{1}(t)|\right) \left(|u_1(t) - u_2(t)| + |\eta_1(t) - \eta_2(t)|\right),\\[1mm]
|w_1(t) - w_2(t)|_W &\le \omega |u_1(t) - u_2(t)| + \gamma |\eta_1(t) - \eta_2(t)|,\\[1mm]
|\dot w_1(t)|_W &\le a(t) + \omega |\dot u_1(t)| + \gamma |\dot \eta_1(t)|\\[1mm]
& \le  \left(1+\frac{\delta}{1-\delta}\right)a(t) + \left(\omega + \frac{\gamma}{1-\delta}\left(1 + \frac{\omega K_1}{c}\right)\right) |\dot u_1(t)|.
\end{align*}
hold for a.\,e. $t\in (0,T)$ By virtue of Hypothesis \ref{gh1}.
We now apply Theorem \ref{2p4} and
insert the above estimates into \eqref{2e11}. The estimate \eqref{eta} now follows with constants $m_0>0$, $m_1>0$ depending only on $c, C_0, C_1, K_0, K_1, \gamma, \omega,C_\xi$, and $C_u$.
\epf

We now prove the main result of this section.

\begin{theorem}\label{3t1}
Let Hypotheses \ref{imh1} and \ref{gh1} hold, and let $u\in \AC$ 
and $x_0 \in X$ be given such that
$x_0 \in Z(g(0,{u(0),} u(0)-x_0))$. Then
there exists a unique solution $\xi \in \Omega(u)$ to \eqref{SGP1}--\eqref{SGP3}.
\end{theorem}

\bpf{Proof}
We proceed by the Banach Contraction Principle. For an arbitrarily given $u \in \AC$ and each $\eta \in\Omega(u)$ we use Lemma \ref{3l2} to find the solution $\xi \in \AC$ of \eqref{GP1}--\eqref{GP3} with $w(t) = g(t,u(t),\eta(t))$. It suffices to prove that the mapping $S:\Omega(u) \to \Omega(u):\eta \mapsto \xi$ is a contraction on $\Omega(u)$.

Let  $\eta_1, \eta_2 \in\Omega(u)$ be given. From \eqref{eta} with
$u_1 = u_2 = u$ it follows that
\begin{align}\label{eta0}
 & |\doxi_1(t) - \doxx(t)| + \frac1c\,
\frac \dd{\dd t} |G(x_1(t),w_1(t)){-}G(x_2(t),w_2(t))| \notag \\ 
& \le \delta |\dot \eta_1(t) - \dot \eta_2(t)|+ 
m(t)\big(|{\xi_1}(t) - {\xi_2}(t)| + |\eta_1(t) - \eta_2(t)|\big)
\end{align}
with $m(t) = m_0(a(t) + b(t) + |\dou(t)|)$. We have $m \in L^1(0,T)$, and we may put for $t \in [0,T]$
$$
M(t) = \int_0^t m(\tau)\dd\tau, \quad M_\ve(t) = \expe^{-\frac1\ve M(t)}, 
$$
where $\ve\in (0,1)$ is chosen in such a way that
\be{eps}
\delta^* := \frac{\delta+\ve}{1-\ve} < 1.
\ee
We multiply \eqref{eta0} by $M_\ve(t)$ and integrate from $0$ to $T$. For simplicity, put
\be{gamma}
\Gamma(t) = \frac1c|G(x_1(t),w_1(t))-G(x_2(t),w_2(t))|. 
\ee
We have
$$
\int_0^T M_\ve(t)\dot \Gamma(t)\dd t = \Gamma(T) M_\ve(T) + \frac1\ve\int_0^T m(t) M_\ve(t)\Gamma(t)\dd t \ge 0, 
$$
hence,
\begin{align}\nonumber
\int_0^T M_\ve(t)|\doxi_1(t) - \doxx(t)|\dd t  
&\le \delta \int_0^T M_\ve(t)|\dot \eta_1(t) - \dot \eta_2(t)|\dd t\\[1mm] \label{eta1}
&\hspace{-20mm} - \ve \int_0^T \dot M_\ve(t) \left(
\int_0^t \big(|{\doxi_1}(\tau) - \doxx(\tau)| + 
|\dot \eta_1(\tau) - \dot\eta_2(\tau)|\big)\dd \tau \right) \dd t.
\end{align}
We use the fact that
$$
M_\ve(T) \left(
\int_0^T \big(|{\doxi_1}(\tau) - \doxx(\tau)|+ 
|\dot \eta_1(\tau) - \dot\eta_2(\tau)|\big)\dd \tau \right) \ge 0,
$$
and integrating by parts in the right-hand side of \eqref{eta1} we obtain
\be{eta2}
(1-\ve)\int_0^T M_\ve(t)|\doxi_1(t) - \doxx(t)|\dd t  \le 
 (\delta + \ve) \int_0^T M_\ve(t)|\dot \eta_1(t) - \dot \eta_2(t)|\dd t.
\ee
Hence, by virtue of \eqref{eps}, the mapping $S: \Omega(u) \to \Omega(u): \eta \mapsto \xi$ is a contraction with respect to the complete metric induced on $\Omega(u)$ by the norm
$$
\|\eta\|_\ve := \int_0^T M_\ve(t)|\dot\eta(t)|\dd t,
$$
and we infer the existence of a solution $\xi \in \Omega(u)$ to \eqref{SGP1}--\eqref{SGP3}.
\epf

\begin{corollary}\label{3c4}
Let Hypotheses \ref{imh1} and \ref{gh1} hold. Then the mapping which with $u\in \AC$ and $x_0 \in Z(g(0,u(0),u(0)-x_0))$ associates the solution $\xi \in \Omega(u)$ of \eqref{SGP1}--\eqref{SGP3} is locally Lipschitz continuous in the sense that for every $R > 2m_0\int_0^T(a(t)+b(t)) \dd t$ there exists $K(R) > 0$ such that if
$u_1, u_2 \in \AC$ are given and
\be{ur}
\int_0^T|\dot u_1(t)|\dd t \le \frac{R}{2m_0} - \int_0^T(a(t)+b(t)) \dd t,
\ee
then the solutions $\xi_i \in \AC$, $i=1,2$ to \eqref{SGP1}--\eqref{SGP3} associated with the inputs $u_i$ and initial conditions $x_0^i\in X$, $x_0^i \in Z(g(0,u_i(0),u_i(0)-x_0^i))$ for $i = 1,2$
satisfy the inequality
\be{3c4-e}
\int_0^T|\doxi_1(t) - \doxx(t)|\dd t \le K(R)  
\left(|{x_0^1} - {x_0^2}| +|{u_1}(0) - {u_2}(0)| + \int_0^T|\dot u_1(t) - \dot u_2(t)|\dd t \right).
\ee
\end{corollary}

\bpf{Proof}
By virtue of \eqref{eta} with $\eta_i = \xi_i$ 
we have for a.\,e. $t\in (0,T)$ that
\begin{align}\nonumber
(1-\delta)|\doxi_1(t) - \doxx(t)| &+ \frac1c\,
\frac \dd{\dd t} |G(x_1(t),w_1(t))-G(x_2(t),w_2(t))| 
\le m_1|\dou_1(t) - \dou_2(t)| \\[1mm]\label{eta3}
&\hspace{-10mm} + 2m_0(a(t) + b(t) + |\dot u_1(t)|) \big(|{u_1}(t) - {u_2}(t)| +
|{\xi_1}(t) - {\xi_2}(t)|\big).
\end{align}
We proceed as in the proof of Theorem \ref{3t1} choosing $\ve \in (0,1-\delta)$ and putting for $t \in [0,T]$ 
\begin{align*}
\hat m(t) &=  2m_0(a(t) + b(t) + |\dot u_1(t)|),\\
\hat M(t) &= \int_0^t \hat m(\tau)\dd\tau,\\
\hat M_\ve(t) &= \expe^{-\frac1\ve \hat M(t)}.
\end{align*}
Note that by \eqref{ur} we have
$$
\hat M(T) = \int_0^T  2m_0(a(t) + b(t) + |\dot u_1(t)|)  \dd t \le R.
$$
Multiplying \eqref{eta3} by $\hat M_\ve(t)$ and using the notation \eqref{gamma} we obtain after integrating from $0$ to $T$ that
\begin{align}\nonumber
&(1-\delta)\int_0^T\expe^{-\frac1\ve \hat M(t)} |\doxi_1(t) - \doxx(t)|\dd t \le 
\Gamma(0) + m_1\int_0^T\expe^{-\frac1\ve \hat M(t)} |\dot u_1(t) - \dot u_2(t)|\dd t 
\\[1mm]
&\quad -\ve \int_0^T \frac{\dd}{\dd t}\hat M_\ve(t)
\left(|{u_1}(0) {-} {u_2}(0)| +
|{\xi_1}(0) {-} {\xi_2}(0)|\right) \dd t \nonumber
\\[1mm] \label{eta4}
&\quad-\ve \int_0^T \frac{\dd}{\dd t}\hat M_\ve(t)\left(\int_0^t|\dot u_1(\tau) {-} \dot u_2(\tau)|+ |\doxi_1(\tau) {-} \doxx(\tau)|\dd \tau \right)\dd t.
\end{align}
On the right-hand side of \eqref{eta4} we integrate by parts and obtain
\begin{align}\nonumber
&(1-\delta-\ve)\int_0^T\expe^{-\frac1\ve \hat M(t)} |\doxi_1(t) - \doxx(t)|\dd t 
\le (m_1+\ve)\int_0^T\expe^{-\frac1\ve \hat M(t)} |\dot u_1(t) - \dot u_2(t)|\dd t
\\[1mm] \label{eta5}
&\qquad + C 
\left(|{u_1}(0) {-} {u_2}(0)| + |{x_0^1} {-} {x_0^2}|\right)
\end{align}
with a constant $C>0$ independent of $R$. Thus, as $0 \le  \hat M(T) \le R$
we obtain the final estimate
\be{eta6}
\int_0^T|\doxi_1(t) - \doxx(t)|\dd t \le C\expe^{R/\ve} \left(|{x_0^1} - {x_0^2}| +|{u_1}(0) - {u_2}(0)| + \int_0^T|\dot u_1(t) - \dot u_2(t)|\dd t \right)
\ee
with a constant $C>0$ independent of $R$, which we wanted to prove.
\epf

\begin{remark} \label{remimp}
The smallness of $\partial_\xi g$ in \eqref{3e1} is indeed a necessary condition for the existence of a solution of the implicit problem even in 1D with $Z(t) = [-r(t), r(t)]$, $r(t) = g(\xi(t))$. It is easy to see that we have existence and uniqueness if $|g'(\xi)| < 1$, and non-existence if $g'(\xi) \le -1$. Furthermore, even in the convex case, the Lipschitz regularity of $\nabla g$ is a necessary condition for uniqueness in the implicit problem. An example of nonuniqueness is provided in \cite[Section 8]{bks} when the
$C^{1,1}$-condition for $g$ is violated.
\end{remark}

\begin{remark} \label{gurson}
Similar result to Theorem \ref{3t1} is obtained if \eqref{im1} is replaced with
\be{imgur}
w(t) = g(t,u(t),V_\xi(t)), \quad V_\xi(t) = \int_0^t |\dot\xi(\tau)\dd\tau.
\ee
In application to elastoplasticity, $V_\xi(t)$ corresponds to dissipated energy during the time interval $[0,t]$, which can be considered as a measure of accumulated cyclic fatigue. For example, in \cite{gurs}, the Gurson model for fatigue is based on the assumption that set $Z(t)$ of admissible stresses shrinks as $V_\xi(t)$ increases.
\end{remark}

\section*{Acknowledgement}
The authors are grateful to Giovanni Colombo for
motivation to study this problem and for stimulating
discussions about the subject. 
They also thank Florent Nacry and Lionel Thibault for pointing out that Proposition \ref{lp1} was already proved in a different way in their paper \cite{AdlHadThi17} with Samir Adly. 
This remark was also made by an anonymous reviewer who referred us to the papers 
\cite{AdlHadThi16, AdlHadThi17}; we are thankful for her/his encouraging report and further comments which helped to improve the manuscript.

\end{document}